\documentclass[onefignum,onetabnum]{siamonline220329}

\newcommand{\abs}[1]{\left| #1 \right|}
\newcommand{\dint}{\mathrm{d}}

\newcommand{\Var}[1]{\mathrm{Var}\left[#1\right]}

\newcommand{\Exp}[1]{\mathbb{E}\left[#1\right]}
\newcommand{\comments}[1]{}

\usepackage{lipsum}
\usepackage{amsfonts}
\usepackage{graphicx}
\usepackage{epstopdf}
\usepackage{algorithmic}
\usepackage{caption}
\usepackage{subcaption}
\usepackage{amsmath,amssymb}

\usepackage{hyperref}

\usepackage{todonotes}

\ifpdf
  \DeclareGraphicsExtensions{.eps,.pdf,.png,.jpg}
\else
  \DeclareGraphicsExtensions{.eps}
\fi

\usepackage{enumitem}
\setlist[enumerate]{leftmargin=.5in}
\setlist[itemize]{leftmargin=.5in}

\newsiamremark{remark}{Remark}
\newsiamremark{hypothesis}{Hypothesis}
\crefname{hypothesis}{Hypothesis}{Hypotheses}
\newsiamthm{claim}{Claim}

\headers{Noise interference to warnings of tipping points}{A. Morr, N. Boers and P. Ashwin}

\title{Internal noise interference to warnings of tipping points\\in generic multi-dimensional dynamical systems}

\author{Andreas Morr\footnote{Corresponding author: \email{andreas.morr@tum.de}} \footnote{Earth System Modelling, School of Engineering and Design, Technical University of Munich, 80333 Munich, Germany} \footnote{Complexity Science, Potsdam Institute for Climate Impact Research, 14473 Potsdam, Germany} \and Niklas Boers\footnotemark[2] \footnotemark[3] \footnote{Department of Mathematics and Statistics, University of Exeter, Exeter EX4 4QF, United Kingdom} \and Peter Ashwin\footnotemark[4]}

\usepackage{amsopn}

\ifpdf
\hypersetup{
  pdftitle={Morr2023NoiseInterference},
  pdfauthor={A. Morr, N. Boers, and P. Ashwin}
}
\fi

\begin{document}

\maketitle
\begin{abstract}
A deterministic dynamical system that slowly passes through a generic fold-type (saddle-node) bifurcation can be reduced to one-dimensional dynamics close to the bifurcation because of the centre manifold theorem.
It is often tacitly assumed that the same is true in the presence of stochasticity or noise so that, for example, critical slowing down (CSD) indicators can be applied as if the system were one-dimensional.
In this work, we show that this is only true when given suitable system observables; specifically, we demonstrate that noise in other dimensions may interfere with indicators of CSD, also referred to as early warning signals (EWS).
We point out a generic mechanism by which both variance and lag-1 autocorrelation [AC(1)], as well as other EWS, can fail to signal an approaching bifurcation.
This can, in principle, occur whenever one noise source drives multiple system components simultaneously.
Even under the favourable assumptions of uncoupled deterministic dynamics and stationary noise, some system observables can then exhibit false negative or false positive CSD indications. We isolate this phenomenon in an example that represents a generic two-dimensional fold-type bifurcation setting.
\end{abstract}



\section{Introduction}\label{sec:Intro}

Many physical systems exhibit abrupt transitions or are suspected to harbour the potential for them \cite{Lenton2012CSDClimate, Ashwin2012tipTypes, Hopcroft2021SaharaAbruptModel, Ren2015CSDElectricity, Boers2022PaleoTipping, Kuehn2015CSDNeuron}.
There has been an increased interest in the analysis and detection of bifurcation-induced transitions \cite{Boers2017DeforAmaz, Boers2021AMOCEWS, Bury2020SpectralEWS, Bury2021CSDMachineLearning, Nazarimehr2020CSDIndicators, Ditlevsen2023AMOCPrediction}.
Such events occur if the system's current equilibrium state destabilises in response to an external change to the underlying dynamics.
The simplest model for a generic bifurcation-induced transition is that of a one-dimensional fold-type bifurcation \cite{Scheffer2009EWSNature, Thompson2010FoldTippingPredict, Guttal2007BistablEcol}.
This approach is reasonable if there is one essential component of the system experiencing negative feedback on its fast time-scale disturbances. 
Moreover, bifurcation theory tells us that this is one of only two possible generic one-parameter bifurcations of an attracting equilibrium, the other being the Hopf bifurcation.
Close enough to a fold-type bifurcation, one can rigorously and generically reduce a multi-dimensional system to a one-dimensional model on a centre manifold.
The negative feedback on this essential component weakens in advance of the critical transition brought about by the external change in dynamics.
The resulting observational characteristic of slower and weaker responses to perturbations is called critical slowing down (CSD).
If the model component undergoing such a bifurcation is a direct observable of the physical system, one can hope to employ the conventional techniques for detecting CSD, i.e.~searching for positive trends in variance or lag-1 autocorrelation [AC(1)] in the time-series data.
These trends are often referred to as early warning signals (EWS) of an impending tipping point in the literature \cite{Kuehn2011CSD, Kuehn2013CSDVar, Lenton2012CSDClimate, Boettner2022CSDCorrNoise}.
These methods are typically motivated by assuming proximity to a fold-type bifurcation of a one-dimensional system whose state is denoted by $x(t)\in\mathbb R$ and whose dynamics depend on a parameter $\alpha$.
In such a case, there will generically be a transformation of the dynamics to a local topological bifurcation normal form that, without loss of generality, is at $x=\alpha=0$ and can be written as:
\begin{equation}\label{eq:1Dtopnorm}
    \dot x= -x^2+\alpha.
\end{equation}
For $\alpha>0$, the system admits one stable equilibrium at $x^*(\alpha)=\sqrt{\alpha}$ and one unstable equilibrium, whereas for $\alpha<0$, it admits no equilibria.

Complex systems that exhibit a clear separation of dynamical time scales can often be effectively represented by a stochastic differential equation (SDE).
This may be motivated, for example, by the projection on macroscopic dynamics \cite{Zwanzig2001Formalism} or by leveraging the ergodic nature of rapidly varying chaotic system components \cite{Melbourne2011ChaosToNoise}.
Adding Gaussian white noise to model the omitted fast dynamics of the original high-dimensional physical system reduces the problem of determining the statistical characteristics (near a stable equilibrium) to the analysis of a one-dimensional Ornstein-Uhlenbeck process.
This is because, in the case of small noise strength $\sigma$, the deterministic dynamics around the equilibrium can be represented by their first-order approximation:
\begin{equation*}
    \dint X_t=-2\sqrt{\alpha}(X_t-x^*(\alpha))\dint t+\sigma\dint W_t,
\end{equation*}
where $W$ is a Wiener process supported on the filtered probability space \linebreak$(\Omega,\mathcal{F},(\mathcal{F}_t)_{t\in\mathbb{R}_+},\mathbb{P})$.
For variance and AC(1) of the observable $X$ sampled at time interval $\Delta t=1$, one finds in the stationary limit
\begin{equation}\label{eq:VarAC1}
    \Var{X}=\frac{\sigma^2}{4\sqrt{\alpha}},\quad\mathrm{AC}_X(1)=\exp(-2\sqrt{\alpha}),
\end{equation}
which both increase as $\alpha>0$ decreases towards the bifurcation at $\alpha=0$.
We may conclude that given time-series observations of a one-dimensional system approaching a fold-type bifurcation, one can reliably detect CSD employing such EWS.
However, the coarse-grained dynamics of real-world systems generally consist of multiple coupled dimensions.
It is {\em a priori} not clear that a given system observable contains information regarding CSD and related EWS of an impending bifurcation.

In a multi-dimensional system, we should expect the noise disturbances within the system to drive multiple components simultaneously, meaning that for a ``typical'' observable being a function of the system state, the noise may interfere with itself and impede any warning of an approach to a tipping point.
In particular, we investigate the implications of multiple independent noise disturbances on critical slowing down in observables near a fold-type bifurcation.

In the following \Cref{sec:centre}, we unpack some of the tacit assumptions underlying the conventional application of EWS methods.
We point out that even though one can isolate the critical component or centre manifold in the multi-dimensional deterministic dynamics, a stochastic coupling of dimensions will persist in general.
This will prove to be problematic when considering non-trivial observables $\Psi$ of the system.
\Cref{sec:OUanalysis} is dedicated to an analytical derivation of the statistical properties of a linear observable $\Psi$ in the two-dimensional case and an assessment of the potential for ``deceitful'' warnings.
Having identified possible pitfalls to the EWS in variance and AC(1), we check that those predicted behaviours indeed manifest in a numerically integrated non-linear system in \Cref{sec:data}.

\section{Centre manifold theory as the motivation for one-dimensional EWS analysis}\label{sec:centre}

The fact that a multi-dimensional system with a generic fold-type bifurcation can be transformed to the topological normal form \cref{eq:1Dtopnorm} is often employed to motivate the use of a single system component in analysing EWS.
To make this statement more concrete and to explore its limitations, we invoke the theory of centre manifolds at a bifurcation point (see \cite[Theorem 5.2]{Kuznetsov2004AppliedBifurcationTheory} and \cite[Theorem 3.4.1]{Guckenheimer1983NonlinearOscillation} for a more detailed discussion).
Consider a dynamical system whose state is denoted by $\mathbf{z}(t)\in\mathbb R^n$, $n\geq 1$, and whose dynamics follow
\begin{equation}\label{eq:nDsysGen}
    \dot{\mathbf{z}}=\mathbf{f}(\mathbf{z},\alpha),
\end{equation}
dependent on a parameter $\alpha\in{\mathbb R}$.
Let us assume that the right-hand side $\mathbf{f}\in C^\infty(\mathbb R^n\times\mathbb R,\mathbb R^n)$ has a bifurcating equilibrium at $(\mathbf{z}_c,\alpha_c)$ that is close to stability.
More precisely, suppose that
\begin{equation}
    \mathbf{f}(\mathbf{z}_c,\alpha_c)=\mathbf 0,\quad \sigma[\nabla_\mathbf{z}\mathbf{f}(\mathbf{z}_c,\alpha_c)]=\{0\} \cup S,
    \label{eq:assump}
\end{equation}
where $\sigma$ denotes the spectrum of the Jacobian $\nabla_{\mathbf{z}} \mathbf{f}$.
Equation \cref{eq:assump} means we assume this spectrum consists of a single zero and of $n-1$ eigenvalues $S\subset {\mathbb C}$ (counting multiplicity) with negative real part. 

Generically, there will be a fold-type bifurcation at this point on varying $\alpha$ \cite{Kuznetsov2004AppliedBifurcationTheory}.
Note that the centre manifold of the system at $(\mathbf{z}_c,\alpha_c)$ extended by a $\dot{\alpha}=0$ is tangent to the two-dimensional plane spanned by the eigenvector associated with the zero eigenvalue of $\nabla_\mathbf{z}\mathbf{f}(\mathbf{z}_c,\alpha_c)$ and the $\alpha$-direction.
This centre manifold is invariant, exponentially attracting and, to linear order, uncoupled from the stable $(n-1)$-dimensional eigenspace.
This implies that the system \cref{eq:nDsysGen} can be reformulated by taking coordinates $\mathbf{z}=(x,\mathbf{y})$ relative to the eigenspaces for $0$ and $S$ and assuming without loss of generality that $(\mathbf{z}_c,\alpha_c)=(\mathbf 0,0)$:
\begin{equation*}\label{eq:locnDsys}
\begin{aligned}
    \dot{x} = &a x^2+b\alpha+\mathcal{O}(2),\\
    \dot{\mathbf{y}} =&C\mathbf{y}+\mathcal{O}(2).
\end{aligned}
\end{equation*}
Here, all eigenvalues of the $(n-1)\times(n-1)$ matrix $C$ have negative real parts, and $a$ and $b$ are normal form coefficients that are generically non-zero.
The $\mathcal{O}(2)$ denotes terms that are quadratic in $x$, $\mathbf{y}$ and $\alpha$. If $a$ and $b$ are non-zero (this corresponds to the genericity conditions (SN2) and (SN3) in \cite[Theorem 3.4.1]{Guckenheimer1983NonlinearOscillation}) then this can be transformed to
\begin{equation}\label{eq:locnDsystrun}
\begin{aligned}
    \dot{x} = &a x^2+b\alpha,\\
    \dot{\mathbf{y}} =&C\mathbf{y}.
\end{aligned}
\end{equation}

With addition of noise, \cref{eq:locnDsystrun} becomes the It\^{o} SDE
\begin{equation}\label{eq:sde-nDsysGen}
    \dint \mathbf{Z}_t=\mathbf{f}(\mathbf{Z}_t,\alpha)\dint t+\Sigma \dint \mathbf{W}.
\end{equation}
Because a Wiener distribution is determined purely by its covariance matrix $\Sigma^T \Sigma $, without loss of generality, and by Cholesky decomposition, we can choose 
\begin{equation*}
    \Sigma= \begin{pmatrix}
    \sigma_x & 0 \\
    \boldsymbol{\sigma}_\mathbf{y}  & \Sigma_\mathbf{y} 
\end{pmatrix}
\end{equation*}
to be a lower triangular matrix.
This gives a system
\begin{equation}\label{eq:sde-locnDsystrun}
\begin{aligned}
    \dint X_t = &[a X_t^2+b\alpha]\dint t + \sigma_x \dint W_t^{(1)},\\
    \dint \mathbf{Y}_t =&[C\mathbf{Y}_t] \dint t + \boldsymbol{\sigma}_\mathbf{y} \dint W_t^{(1)}+\Sigma_\mathbf{y}  \dint \mathbf{W}_t^{(2)},
\end{aligned}
\end{equation}
where $\sigma_{x}$ is a scalar, $\boldsymbol{\sigma}_\mathbf{y} $ is a vector, $\Sigma_\mathbf{y} $ is a lower triangular $(n-1)\times(n-1)$ matrix, and $W^{(1)}$ and $\mathbf{W}^{(2)}$ are standard Wiener processes that are independent in each component.
Note that we choose here a lower triangular $\Sigma$ so that $\dint X$ is formulated with a single white-noise term.
Later on, we will equivalently choose an upper triangular $\Sigma$ for analysis. 
A more general noise model would include a state dependence of $\Sigma$.
In this work, we restrict ourselves to the case of an approximately constant $\Sigma$.
The application of the above normal-form theory demands low-amplitude noise.

An important implication of Eq.~\cref{eq:sde-locnDsystrun} is that, close to the bifurcation, the $X$ coordinate functions as an optimal observable, i.e., it exhibits the expected behaviour with respect to CSD and its expression of EWS approximated in Eq.~\cref{eq:VarAC1}.
This confirms it is reasonable to expect that there is a destabilising system component which exhibits EWS: \cite{Held2004CSDDetection, Dakos2018MultDResiliencevsEigenvalues_EOF, Ghadami2020FindLowestEigEWS, VanNes2007SlowRecoveryCSD, George2023MultiVarEWSOverview}.
The reduction of the deterministic dynamics seems to offer a valuable simplification: The critical dimension or centre manifold dynamics in $X$ are uncoupled from all other coarse-grained equilibrium dynamics subsumed in $\mathbf Y$ under such a transformation.
However, this simplification is deceptive, as it does not consider possible couplings in the stochastic fast dynamics.
We explore this generic circumstance in more detail for a two-dimensional system.

Consider the case $n=2$ (so $y\in{\mathbb R}$) 
and assume $a=-1$, $b=1$ and $C=-1$ so that:
\begin{equation*}\label{eq:2Dsys}
\begin{aligned}
    \dot x=&-x^2+\alpha,\\
    \dot y=&-y.
\end{aligned}
\end{equation*}
Adding a generic white-noise term yields the following SDE with a non-diagonal (anisotropic) noise coupling matrix:
\begin{equation}\label{eq:2DsysStoch}
    \begin{pmatrix}\dint X_t\\\dint Y_t\end{pmatrix}=\begin{pmatrix}-X_t^2+\alpha\\-Y_t\end{pmatrix}\dint t+\Sigma\begin{pmatrix}\dint W^{(1)}_t\\\dint W^{(2)}_t\end{pmatrix}.
\end{equation}

One can arrive at Eq.~\cref{eq:2DsysStoch} by transforming some original model to show a centre manifold which is isolated from the remaining deterministic dynamics, as we have discussed above.
No such isolation can be simultaneously achieved on the side of the stochastic dynamics.
One can hope to observe CSD in the $X$ component of the system, as this sub-system is equivalent to the one-dimensional case discussed in \Cref{sec:Intro}.
However, in applications, one may not be free to choose a system observable.
An available observable will be of the form $\Psi_t=g(X_t,Y_t)$, where $g$ could be any non-linear function.
We show that already the subset of simple linear observables
\begin{equation}\label{eq:Psi}
    \Psi_t:=\cos(\beta)X_t+\sin(\beta)Y_t
\end{equation}
proves problematic in the context of anticipating critical transitions using EWS.
Equation \cref{eq:Psi} covers all possible linear combinations of $X_t$ and $Y_t$, and the phenomenon we discuss is impervious to an absolute scaling of the observable.

We investigate the linearised system corresponding to Eq~\cref{eq:2DsysStoch} around $(x^*(\alpha),y^*(\alpha))=(\sqrt{\alpha},0)$.
If the noise amplitude is small, this will be an exceptionally good approximation of the full non-linear dynamics, as the system remains largely in the linear regime:
\begin{equation}\label{eq:2DsysLin}
    \begin{pmatrix}\dint X_t\\\dint Y_t\end{pmatrix}=\begin{pmatrix}-2\sqrt{\alpha}&0\\0&-1\end{pmatrix}\begin{pmatrix}X_t-x^*(\alpha)\\Y_t\end{pmatrix}\dint t+\begin{pmatrix}\sigma_X&c\\0&\sigma_Y\end{pmatrix}\begin{pmatrix}\dint W^{(1)}_t\\\dint W^{(2)}_t\end{pmatrix}.
\end{equation}
When defining the noise term through the matrix $\Sigma$, we are given three degrees of freedom in the form of the respective noise couplings $\sigma_X$, $\sigma_Y$ and $c$. 
Note that any matrix with an additional coupling of the white-noise term $\dint W^{(1)}$ to $\dint Y$ can always be renormalised to an upper-triangular matrix.
This is equivalent to the lower-triangular formulation of the noise term in Eq.~\cref{eq:sde-locnDsystrun}.
Similarly, any third white-noise term could always be represented in the above model equation containing two white-noise terms.
We call $c$ the cross-coupling of the noise, which, as discussed above, will be non-zero in the general case.
We will see that introducing only this additional coupling to the model can have grave implications on the ability to detect the approach of a bifurcation through conventional CSD methods.

\section{Statistical properties of the observable $\Psi$ under linearised dynamics}\label{sec:OUanalysis}
To generalise the linear setting to applications beyond the particular parametrisation of the normal form given in the previous section, we now write the 2D-Ornstein-Uhlenbeck process \cref{eq:2DsysLin} in question generally as
\begin{align}
    \begin{pmatrix}\dint X_t\\\dint Y_t\end{pmatrix}&=\begin{pmatrix}-\lambda_X&0\\0&-\lambda_Y\end{pmatrix}\begin{pmatrix}X_t\\Y_t\end{pmatrix}\dint t+\begin{pmatrix}\sigma_X&c\\0&\sigma_Y\end{pmatrix}\begin{pmatrix}\dint W^{(1)}_t\\\dint W^{(2)}_t\end{pmatrix}\label{eq:OU}\\
    &=:A\begin{pmatrix}X_t\\Y_t\end{pmatrix}\dint t+\Sigma\,\begin{pmatrix}\dint W^{(1)}_t\\\dint W^{(2)}_t\end{pmatrix}\nonumber
\end{align}
The deterministic linear dynamics are still uncoupled, and we have carried along the cross-coupling $c$ of the noise terms.
Without loss of generality, we have also centred the dynamics of the $X$ dimension around $0$.
In this more abstract formulation, we are interested in the statistical properties of the observable $\Psi$ and their dependence on the value of $\lambda_X$.
This is because, while approaching the fold-type bifurcation, $\lambda_X$ will generically decrease towards zero and canonically incur CSD in the $X$ dimension.
Whether this CSD is also detectable in the observable $\Psi$ is the main subject of this work.
The two-dimensional SDE \cref{eq:OU} together with the initial condition $(X_0,Y_0)^\top=(0,0)^\top$ has the following closed-form solution \cite{Gardiner1985StochAna}:
\begin{equation*}
    \begin{pmatrix}X_t\\Y_t\end{pmatrix}=\int_0^t\exp((t-s)A)\Sigma\begin{pmatrix}\dint W_s^{(1)}\\\dint W_s^{(2)}\end{pmatrix}.
\end{equation*}
Consider the covariance matrix $V$ of this zero-mean process, i.e.
\begin{align*}
    V:&=\Exp{\left(X_t,Y_t\right)^\top\left(X_t,Y_t\right)}\\&=\Exp{\int_0^t\int_0^t\exp((t-s_1)A)\Sigma\begin{pmatrix}\dint W_{s_1}^{(1)}\\\dint W_{s_1}^{(2)}\end{pmatrix}\left(\dint W_{s_2}^{(1)},\dint W_{s_2}^{(2)}\right)\Sigma^\top\exp((t-s_2)A)^\top}\\
    &=\int_0^t\exp((t-s)A)\Sigma\Sigma^\top\exp((t-s)A)^\top\dint s,
\end{align*}
where we have employed the It\^o isometry in the last step.
Taking the time derivative on both sides and investigating the stationary limit yields
\begin{equation*}
    \frac{\dint}{\dint t}V=AV+VA^\top+\Sigma\Sigma^\top=0.
\end{equation*}
For the system in \eqref{eq:OU}, this so-called Lyapunov matrix-equation \cite{Barnett1968LyapunovMatrixEqApplications} is solved by
\begin{equation*}
    V=\begin{pmatrix}\frac{\sigma_X^2+c^2}{2\lambda_X}&\frac{\sigma_Yc}{\lambda_X+\lambda_Y}\\\frac{\sigma_Yc}{\lambda_X+\lambda_Y}&\frac{\sigma_Y^2}{2\lambda_Y}\end{pmatrix}.
\end{equation*}
Adopting the initial condition $(X_0,Y_0)^\top\sim\mathcal{N}(0,V)$, the process is stationary, and the time-covariance matrix for all $\tau\geq0$ can be similarly deduced as
\begin{equation*}
    R(\tau):=\Exp{\left(X_t,Y_t\right)^\top\left(X_{t+\tau},Y_{t+\tau}\right)}=\exp(\tau A)V.
\end{equation*}
Note that $R(\tau)$ for $\tau\neq0$ will in general not be a symmetric matrix, since $\Exp{X_tY_{t+\tau}}=\Exp{X_{t+\tau}Y_t}$ if and only if $c=0$ or $\lambda_X=\lambda_Y$:
\begin{equation*}
    R(\tau)=\begin{pmatrix}\frac{\sigma_X^2+c^2}{2\lambda_X}\exp(-\lambda_X\abs{\tau})&\frac{\sigma_Yc}{\lambda_X+\lambda_Y}\exp(-\lambda_X\abs{\tau})\\\frac{\sigma_Yc}{\lambda_X+\lambda_Y}\exp(-\lambda_Y\abs{\tau})&\frac{\sigma_Y^2}{2\lambda_Y}\exp(-\lambda_Y\abs{\tau})\end{pmatrix} \,.
\end{equation*}
The same holds true for the time-correlation matrix with entries $r(\tau)_{i,j}:=R(\tau)_{i,j}/\sqrt{V_{i,i}V_{j,j}}$:
\begin{equation*}
    r(\tau)=\begin{pmatrix}\exp(-\lambda_X\abs{\tau})&\frac{2c}{\sqrt{\sigma_X^2+c^2}}\frac{\sqrt{\lambda_X\lambda_Y}}{\lambda_X+\lambda_Y}\exp(-\lambda_X\abs{\tau})\\\frac{2c}{\sqrt{\sigma_X^2+c^2}}\frac{\sqrt{\lambda_X\lambda_Y}}{\lambda_X+\lambda_Y}\exp(-\lambda_Y\abs{\tau})&\exp(-\lambda_Y\abs{\tau})\end{pmatrix}.
\end{equation*}
For any linear observable $\Psi(\beta)=\cos(\beta)X+\sin(\beta)Y$ of the 2D-system, we may now compute the variance and the autocorrelation function in the stationary limit:
\begin{equation*}
    \Var{\Psi}=\cos^2(\beta)V_{1,1}+\sin^2(\beta)V_{2,2}+2\cos(\beta)\sin(\beta)V_{1,2},
\end{equation*}
where $V_{i,j}$ is the respective entry of the covariance matrix $V$.
For the autocorrelation in time, we have:
\begin{align*}
    \mathrm{AC}_\Psi(1)&:=\Exp{\Psi_t\Psi_{t+\tau}}/\Var{\Psi}\\=(\cos^2&(\beta)R(\tau)_{1,1}+\sin^2(\beta)R(\tau)_{2,2}+\cos(\beta)\sin(\beta)(R(\tau)_{1,2}+R(\tau)_{2,1}))/\Var{\Psi}.
\end{align*}
As expected from the uncoupled nature of the deterministic dynamics, if we choose $\beta=0$, i.e. $\Psi=X$, we essentially retrieve the conventional quantities from the one-dimensional case.
This observable should be considered as the optimal observable of the system.
Also, not surprisingly, if we only inspect the second dimension $\Psi=Y$, the variance and AC(1) will not depend on the destabilisation of $X$, i.e. $\lambda_X\rightarrow0$, and can not function as EWS.
See \cref{tab:edgeCases} for explicit expressions.

\begin{table}[H]
\footnotesize
  \caption{Edge cases of the linear observable setting. The conventional EWS information is perfectly conserved in the first case and entirely lost in the second.}\label{tab:edgeCases}
\begin{center}
  \begin{tabular}{|c|c|c|c|} \hline
    $\beta$ & \bf$\Psi$ & \bf$\Var{\Psi}$ & \bf$\mathrm{AC}_\Psi(1)$              \\ \hline
0       & $X$    & $(\sigma_X^2+c^2)/2\lambda_X$          & $\exp(-\lambda_X)$ \\
$\pi/2$ & $Y$    & $\sigma_Y^2/2\lambda_Y$              & $\exp(-\lambda_Y)$ \\ \hline
  \end{tabular}
\end{center}
\end{table}

In all other cases, the observable $\Psi$ consists of some amplitude of $X$ and $Y$, i.e. $\beta\neq k\pi/2$ for $k\in\mathbb{Z}$.
Even here, the EWS seem to be reliable at first glance. We have
\begin{align*}
    \Var{\Psi}\xrightarrow{\lambda_X\rightarrow0}\infty,\quad\mathrm{AC}_\Psi(1)\xrightarrow{\lambda_X\rightarrow0}1.
\end{align*}
These increases are monotonic in proximity to the bifurcation, where $\lambda_X=0$.
More precisely, for any parameter setting $\lambda_Y>0$ and $\sigma_X,\sigma_Y,c\neq 0$, there exists a $\lambda_X^*>0$ with the property that both variance and AC(1) are monotonically increasing as $\lambda_X$ decreases from $\lambda_X^*$ to $0$.
In that sense, there still is a time span in the advent of the fold-type bifurcation, in which the linearised dynamics predict a clear increase in both variance and AC(1), making them applicable as EWS in this case.

\begin{figure}[H]
\centering
\includegraphics[width=\textwidth]{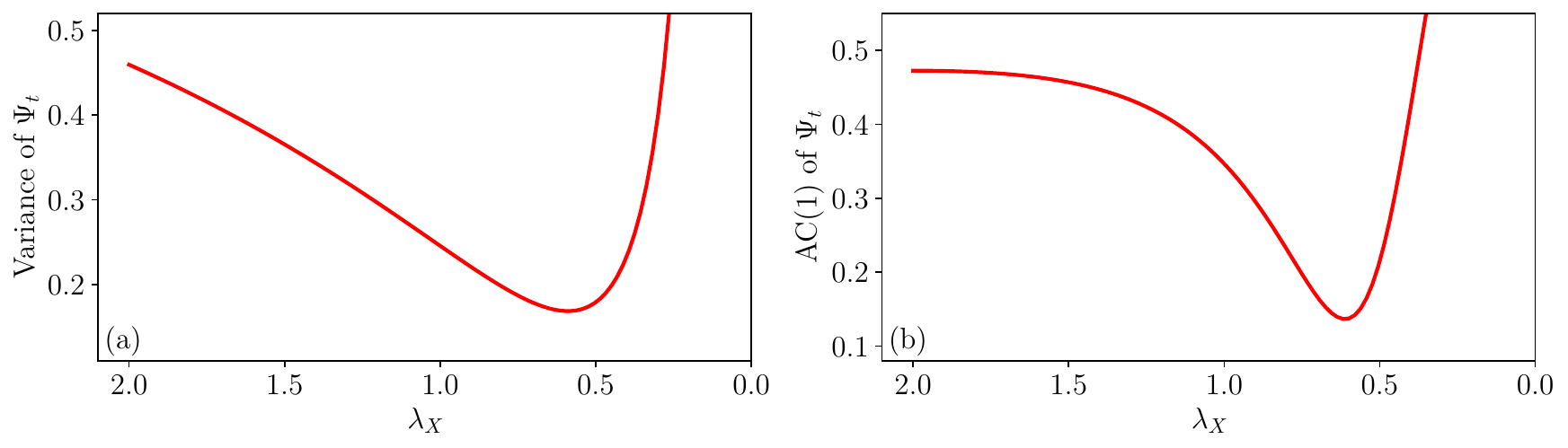}
\caption{\textbf{(a)} Variance and \textbf{(b)} AC(1) of the observable $\Psi$ for $\beta=-\frac{\pi}{4}$ in the linearised setting of Eq.~\cref{eq:exOU} described in the main text.
While these two classical CSD indicators ultimately increase to their expected limits $\infty$ and $1$ respectively as $\lambda_X\rightarrow0$, the behaviour up until that increase can be deceitful with respect to critical slowing down.}
\label{fig:varac1}
\end{figure}

However, the practical caveat to this finding is the value of $\lambda_X^*$ and the behaviour of the EWS before this threshold.
We give the following example, which illustrates that using EWS has potential pitfalls and may lead to a misinterpretation of the actual destabilisation in the underlying dynamics.
For this, choose $\lambda_Y=1$, $\sigma_X=0.1$, $\sigma_Y=2$ and $c=1$, meaning
\begin{align}
    \begin{pmatrix}\dint X_t\\\dint Y_t\end{pmatrix}&=\begin{pmatrix}-\lambda_X&0\\0&-1\end{pmatrix}\begin{pmatrix}X_t\\Y_t\end{pmatrix}\dint t+\begin{pmatrix}0.1&1\\0&2\end{pmatrix}\begin{pmatrix}\dint W^{(1)}_t\\\dint W^{(2)}_t\end{pmatrix},\label{eq:exOU}
\end{align}
and examine the observable $\Psi(\beta)$ associated with the mixing parameter $\beta=-\frac{\pi}{4}$, i.e. \linebreak$\Psi_t=2^{-1/2}(X_t-Y_t)$.
\cref{fig:varac1} shows the variance and AC(1) of $\Psi$ as a function of the destabilising eigenvalue $\lambda_X$.
Although the expected increasing behaviour can be observed when decreasing $\lambda_X$ from $\lambda_X^*\approx0.6$ towards $0$, both CSD indicators, in fact, decrease as $\lambda_X$ decreases toward $\lambda_X^*$.
It is this veiling of the critical slowing down of the $X$ dimension that we would like to point out in this work.

The underlying prerequisite of the phenomenon explored above is the following configuration: The noise component $\dint W^{(2)}$ is positively coupled to both dimensions of the deterministic dynamics.
At the same time, the observable $\Psi$ is defined such that the two linear components have opposing signs.
This leads to the disturbances in $X$ and $Y$ interfering with each other in the summary observable $\Psi$.
The appendix gives a result stating that this effect is ubiquitous with respect to the specific parameter choices.
Theorem~\ref{thm:Theorem} implies that, as long as we have a cross-coupling $c\neq0$, the phenomenon persists for all generic choices of $\sigma_X$, $\sigma_Y$, $c$ and $\lambda_Y$, and it can occur in arbitrarily small proximity to the bifurcation point.
On the other hand, the extent to which it quantitatively inhibits the use of EWS in an application setting will heavily depend on these values.

\begin{figure}[H]
\centering
\includegraphics[width=\textwidth]{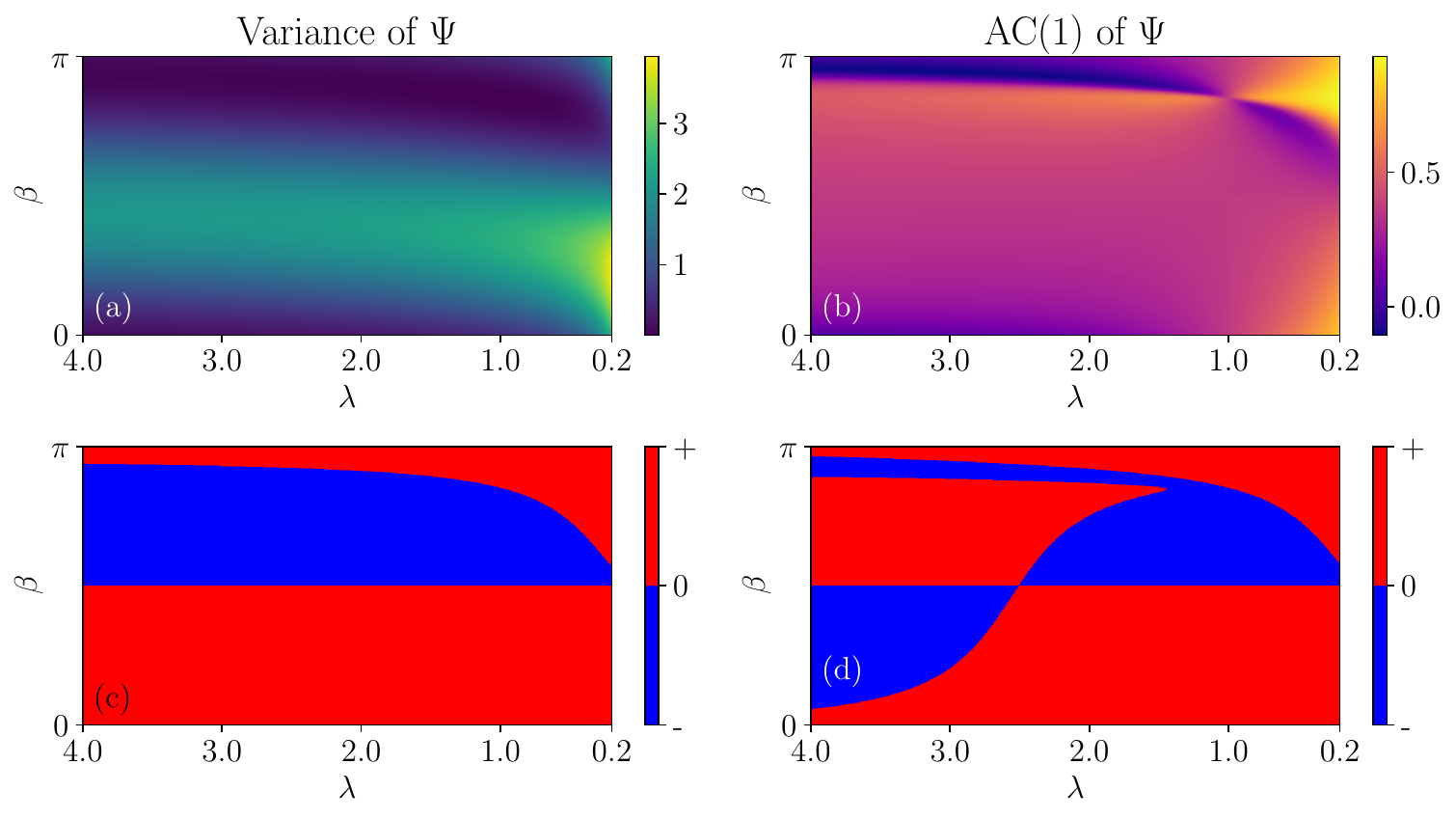}
\caption{\textbf{(a)} variance and \textbf{(b)} AC(1) of $\Psi$ for different values of $\beta$ and $\lambda_X$ in the 2D-Ornstein-Uhlenbeck setting of Eq.~\cref{eq:exOU}. For $\beta=\frac{3\pi}{4}$ one obtains \cref{fig:varac1}, since variance and AC(1) of $\Psi$ are unchanged by a transformation $\beta+k\pi$, $k\in\mathbb Z$.
In \textbf{(c)} and \textbf{(d)}, the colour signifies the trend of variance and AC(1), respectively, as one decreases $\lambda_X$ from left to right during CSD. The plot is red, where this trend is positive, and blue, where it is negative.
In this sense, all blue regions represent the values of $\beta$ and $\lambda_X$ where deceitful CSD indications are prevalent.}
\label{fig:planar}
\end{figure}

Giving a general guideline of which configurations exhibit a lengthy and pronounced decrease in variance and AC(1), even though CSD would imply an increase, is not feasible given the number of degrees of freedom.
We provide a Python script with which the prevalence of the effect can be assessed for different parametrisations of the discussed two-dimensional model.
Central to this script is the plotting of the approach to the critical value $\lambda_X^{\text{crit}}=0$ for different values of $\beta$ and $\lambda_X$(see \cref{fig:planar} and its description).

\section{Observations in model data}\label{sec:data}
\begin{figure}[t]
\centering
\includegraphics[width=\textwidth]{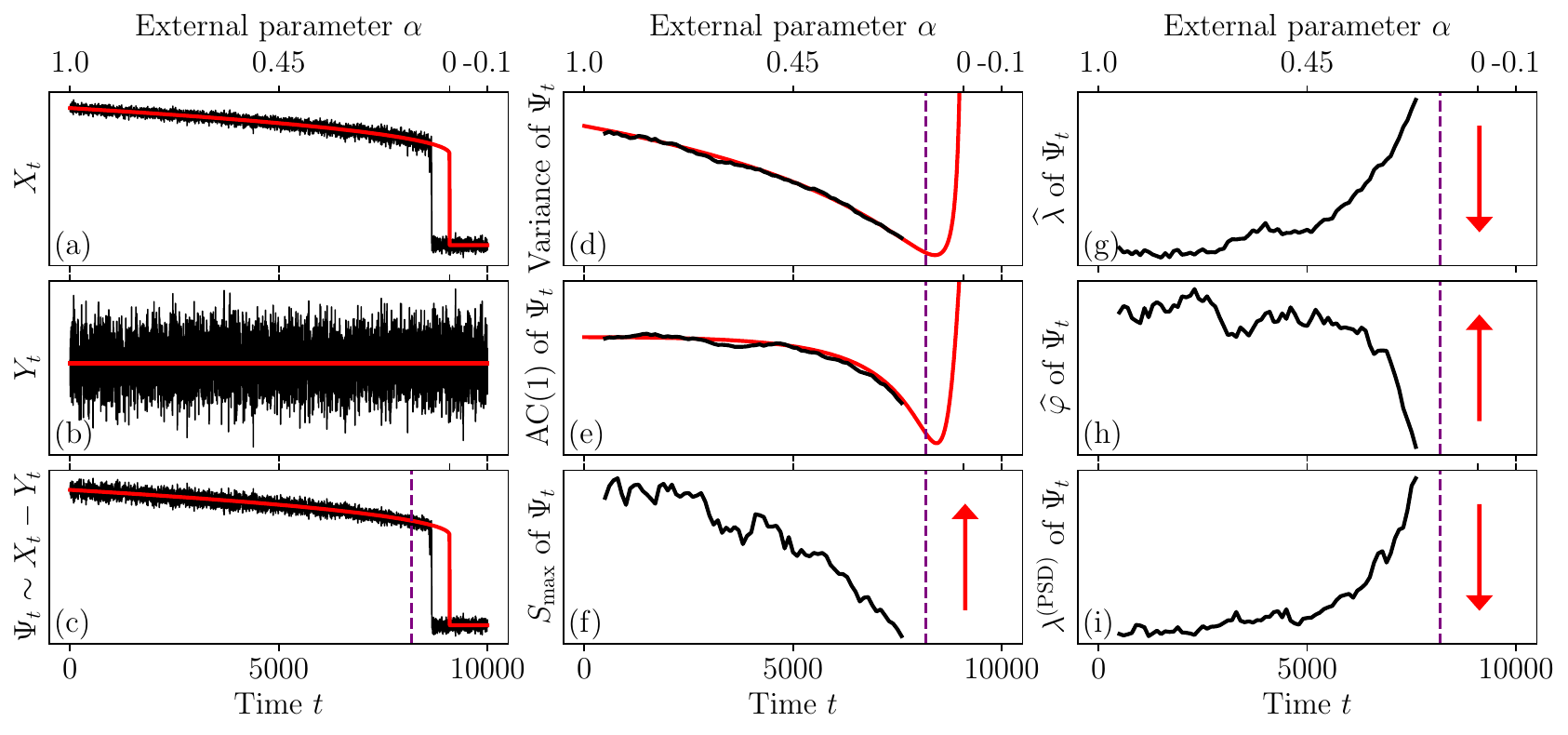}
\caption{\textbf{(a)} - \textbf{(c)} Sample paths of $X_t$, $Y_t$ and $\Psi_t$ following the SDE \cref{eq:foldSDE} in the main text with $\beta=-\frac{\pi}{4}$.
Up until the purple dashed line, after which noise-induced tipping may frequently come into play, the variance and AC(1) of $\Psi$ are each estimated over windows of length $N=1000$. The averaged results over 20 runs are plotted in \textbf{(d)} and \textbf{(e)}, respectively, along with the theoretically computed values from the linearised model. Their decrease would be misinterpreted as a stabilisation of the system rather than an approach to an abrupt transition. In \textbf{(f)} - \textbf{(i)}, alternative indicators for CSD are calculated. These are \textbf{(f)} the maximum of the observed power spectral density (PSD) \cite{Bury2020SpectralEWS}, \textbf{(g)} the linear restoring rate estimated through a Langevin approach \cite{Morr2024KramersMoyalEWS}, \textbf{(h)} a generalised least squares estimator for an AR(1) model driven by correlated noise \cite{Boers2021AMOCEWS} and \textbf{(i)} the linearised restoring rate estimated via the PSD \cite{Morr2024RedNoiseCSD}. These quantities, too, exhibit the respective opposite trend of what would be expected of a destabilising system (indicated by the red arrows).}
\label{fig:foldobservable}
\end{figure}

After identifying potential pitfalls to employing variance and AC(1) on linear observables of the linearised system \cref{eq:2DsysLin}, we test these predictions on synthetic data of the true fold-type bifurcation in Eq.~\cref{eq:2DsysStoch}.
To this end, we integrate the following system of SDEs using the Euler-Mayurama scheme:
\begin{align}\label{eq:foldSDE}
\begin{split}
    \begin{pmatrix}\dint X_t\\\dint Y_t\end{pmatrix}&=\begin{pmatrix}-X_t^2+\alpha(t)\\-Y_t\end{pmatrix}\dint t+\varepsilon\begin{pmatrix}0.1&1\\0&2\end{pmatrix}\begin{pmatrix}\dint W^{(1)}_t\\\dint W^{(2)}_t\end{pmatrix},\\
    \alpha(t)&=1-\frac{11}{10}\frac{t}{T}.  
\end{split}
\end{align}
Over the time span of $T=10^4$, this integration is performed at a time-step of $\delta t=1/30$ and subsequently sampled at every 30th value.
In order to ensure that premature noise-induced tipping is restricted to a relatively short period before the bifurcation, we introduce the noise-scaling parameter $\varepsilon=0.1$.
For better illustration, we artificially introduce a second stable state for $X$ past the bifurcation point.
These dynamics are omitted in the SDE \cref{eq:foldSDE}.
The results for the linear observable
\begin{equation*}
    \Psi_t=\cos(\beta)X_t+\sin(\beta)Y_t=\frac{1}{\sqrt{2}}(X_t-Y_t)
\end{equation*}
with $\beta=-\frac{\pi}{4}$ are given in \cref{fig:foldobservable}.
Even though the dynamics are now non-linear and non-autonomous, the contemporary linearised system still allows for a good assessment of the expected time-series statistics.
The corresponding estimators of variance and AC(1) indeed give the expected results.
This is due to the relatively small noise strength, which keeps the system close to the equilibrium point and thus in an approximately linear regime.
The non-autonomous nature of the system is negligible, as the rate of change introduced in $\alpha(t)$ is slow in comparison with the internal variability induced by the driving noise and the time scales of the restoring dynamics.
The investigated time series windows are short enough with respect to the change in $\alpha$, such that the true variance and AC(1) are approximately constant during that time interval.
The quantities decrease with time, even though the destabilisation of the system in the dimension of the $X$ component should conventionally incur an increase.
This is caused by an interference of the shared driving noise in the observable $\Psi$.
Note that similarly, a stabilisation of the system, i.e.~an increase in $\lambda_X$, would cause the variance and AC(1) to increase.
This is a setting in which a false alarm with respect to CSD would be raised.

\section{Discussion}
We consider a generic fold-type bifurcation in a class of multi-dimensional dynamical systems, where we assume that noise drives multiple components of the system.
Even though the centre manifold theorem yields one distinct, deterministically uncoupled dimension exhibiting pure CSD, this optimal observable will generally not be available for analysis.
We demonstrate that such a setup can easily produce observables that, because of internal noise interference, exhibit deceiving trends in variance and lag-1 autocorrelation with respect to CSD.
False or missed alarms in EWS have been shown to be an issue in systems driven by non-stationary or non-white noise \cite{Dutta2018RobustnessEWScolored, Kuehn2022ColourBlind, Boers2021AMOCEWS, Boettner2022CSDCorrNoise, Proverbio2023OptEWSManyRefs, Morr2024RedNoiseCSD, Morr2024KramersMoyalEWS} or systems exhibiting deterministic coupling \cite{Bury2020SpectralEWS}.
It has also been shown that a misrepresentation of the natural system as a fold-type bifurcation with direct noise influence \cite{Boerlijst2013SilentCollapse, Kefi2013EWSforNonCatastrophic, Dakos2012VarianceCanDecrease} or issues with the availability of time series data \cite{Ditlevsen2010SignalNoise, Benyami2023AMOCDataCSD, Smith2023ReliabilityResilienceEstimation} may lead to similar problems.
This work introduces a previously unknown complication with the detection of CSD via variance and AC(1), as well as several other techniques introduced in Refs.~\cite{Bury2020SpectralEWS, Morr2024KramersMoyalEWS, Boers2021AMOCEWS, Morr2024RedNoiseCSD, Clarke2023ROSA} (see also \cref{fig:foldobservable}f-i).
Even though the severity of the effect described here depends heavily on the relative values of the system's restoring rates and noise couplings, its presence is ubiquitous (see \Cref{thm:Theorem}).
We have also shown that in close proximity to the bifurcation point, EWS can still be expected to deliver information about the impending transition.
However, since non-zero noise always incurs the possibility of noise-induced tipping well before the actual vanishing of a contemporary equilibrium, the discussed complication poses a serious challenge to the reliability of EWS \cite{Ashwin2012tipTypes}.
Furthermore, if present, the effect can only be ruled out by choice of a different observable $\Psi$, which is often not achievable in the observations of natural systems.

Whenever possible, a sufficient physical understanding and previous model validation should confirm that the conventional assumption of a directly observed noisy dynamical system along its centre manifold is warranted.
In this case, the presented results do not call into question the expressiveness of conventional EWS.

\section*{Data Availability}
Visit the GitHub repository \href{https://github.com/andreasmorr/LinearObservablesCSD.git}{LinearObservablesCSD} to access the code used to generate all figures in this manuscript.

\section*{Acknowledgments}
This project has received funding from the European Union's Horizon Europe research and innovation programme under grant agreement No. 101137601 (ClimTip). This is ClimTip contribution \#14. N.B. acknowledges further funding from the Volkswagen Stiftung and the European Union’s Horizon 2020 research and innovation programme under the Marie Sklodowska-Curie grant agreement No. 956170. Funded by the European Union. Views and opinions expressed are however those of the author(s) only and do not necessarily reflect those of the European Union or the European Climate, Infrastructure and Environment Executive Agency (CINEA). Neither the European Union nor the granting authority can be held responsible for them.  For the purpose of open access, the authors have applied a Creative Commons Attribution (CC BY) licence to any Author Accepted Manuscript version arising from this submission.

\bibliographystyle{siamplain}
\bibliography{BibAll}

\appendix

\section{Noise interference of early warning signals in a 2D linear SDE}

For completeness, we state and prove the following result about the ubiquity of noise interference to early warning signs.

\begin{theorem}\label{thm:Theorem}
    Consider the linear SDE \cref{eq:OU} with non-zero constants $\lambda_Y>0$, $\sigma_X$, $\sigma_Y$ and $c$. Choose any small $\delta>0$ and consider the decrease of $\lambda_X$ from $\lambda_Y$ down to $\delta$, representing a destabilisation of the system. There exists an open interval of mixing parameters $\beta\in(\beta^1,\beta^2)$ for which both variance and AC(1) of $\Psi$ defined through Eq.~\cref{eq:Psi} decrease with decreasing $\lambda_X$.
\end{theorem}
\indent\textit{Proof.} We can find such an interval of mixing parameters $\beta\in (\beta^1,\beta^2)$ close to the choice $\Psi=Y$, i.e. $\beta=\pi/2$. We have already derived analytical expressions of $\mathrm{Var}[\Psi]$ and $\mathrm{AC}_{\Psi}(1)$. They are smooth functions of the system's parameters. The quantities decrease with a decreasing $\lambda_X$ if their partial derivative with respect to $\lambda_X$ is positive. We find at the trivial choice of $\Psi=Y$ that for all $\lambda_X\in\mathbb R$
\begin{equation}\label{eq:deriv}
    \partial_{\lambda_X}\mathrm{Var}[\Psi]\Big|_{\beta=\pi/2}=0 \quad\mathrm{and}\quad \partial_{\lambda_X}\mathrm{AC}_{\Psi}(1)\Big|_{\beta=\pi/2}=0.
\end{equation}
We find that the derivative of these quantities with respect to the mixing parameter $\beta$ is non-zero and equal in sign for both of the quantities:
\begin{align*}
    \partial_{\beta}\partial_{\lambda_X}\mathrm{Var}[\Psi]\Big|_{\beta=\pi/2}&=\frac{-2c\sigma_Y}{(\lambda_X+\lambda_Y)^2}\\
    \partial_{\beta}\partial_{\lambda_X}\mathrm{AC}_{\Psi}(1)\Big|_{\beta=\pi/2}&=\frac{-2c\lambda_Y\left(e^{-\lambda_X}\left(1+\lambda_X+\lambda_Y\right)-e^{-\lambda_Y}\right)}{(\lambda_X+\lambda_Y)^2\sigma_Y}
\end{align*}
This sign is equal to $-1\cdot \mathrm{sgn}(c\sigma_Y)$, and it is uniform across all choices of $\lambda_X\leq\lambda_Y$.
Consider without loss of generality a positive sign.
Then for the compact interval $\lambda_X\in[\delta,\lambda_Y]$ we can find $\gamma_\delta>0$ such that the quantities in Eq.~\cref{eq:deriv} are positive
on the interval $(\pi/2,\pi/2+\gamma_{\delta})$.
The counterfactual to this statement implies the existence of a value $\widehat\lambda_X$ with, e.g.,
\begin{equation*}
    \partial_{\lambda_X}\mathrm{Var}[\Psi]\Big|_{\beta=\pi/2,\,\lambda_X=\widehat\lambda_X}=0 \quad\mathrm{and}\quad \partial_{\beta}\partial_{\lambda_X}\mathrm{Var}[\Psi]\Big|_{\beta=\pi/2,\,\lambda_X=\widehat\lambda_X}>0
\end{equation*}
while
\begin{equation*}
    \max_{\beta\in[\pi/2,\gamma]}\left\{\partial_{\lambda_X}\mathrm{Var}[\Psi]\Big|_{\beta=\pi/2,\,\lambda_X=\widehat\lambda_X}\right\}=0 \text{ for some }\gamma>0,
\end{equation*}
a contradiction.
For a negative sign, an analogous interval $(\pi/2-\gamma_{\delta},\pi/2)$ can be found.

This interval constitutes the desired range of mixing parameters.
Note that the chosen interval is not the only region for which the phenomenon discussed in this work is prevalent (see e.g.~\cref{fig:planar}c and d).
\hfill$\square$

\end{document}